\documentclass[12pt]{amsart}

\usepackage{amssymb, amsmath}
\usepackage{amsthm}
\usepackage{graphicx}
\usepackage{hyperref}
\hypersetup{
    colorlinks=true,
    linkcolor=blue,
    urlcolor=blue,
    }
\usepackage{comment}
\usepackage{amsrefs}

\usepackage{caption}
\usepackage{subcaption}
\usepackage{xcolor}

\theoremstyle{plain} 
\newtheorem{thm}{Theorem}[section]
\newtheorem{lem}{Lemma}[section]
\newtheorem{prop}{Proposition}[section]

\theoremstyle{definition}
\newtheorem{dfn}[thm]{Definition}

\theoremstyle{remark}
\newtheorem{rmk}[thm]{Remark}

\newcommand{\CP}{\mathbf C P^2}
\newcommand{\Mod}{\text{Mod}}

\begin{document}
\title[Genus-2 factorizations of type (4,3)]{Equivalent genus-2 factorizations of type $(4,3)$}

\author[E. Huang]{Evan Huang}
\address{Department of Mathematics\\
Princeton University\\
 Princeton, NJ, 08544}
\email{evan.huang@princeton.edu}

\maketitle

\section{Introduction}

In \cite{oz}, Ozbagci showed that the minimal number of singular fibers in a genus-2 Lefschetz fibration over $S^2$ is either 7 or 8. Xiao \cite{xiao} gave a construction coming from algebraic geometry that shows the minimal number is 7, although he did not write down the vanishing cycles of his fibration. Birman-Hilden \cite{birman} showed that $\Mod_2$, the mapping class group of a genus-2 surface, has abelianization $\mathbf Z/10\mathbf Z$. This implies that any positive factorization of length 7 in $\Mod_2$ has 4 nonseparating curves and 3 separating curves (see \cite{bk}, Lemma 5). Such a factorization is said to have type $(4,3)$.

In \cite{bk}, Baykur-Korkmaz constructed a positive factorization of length 7 in $\Mod_2$ using a combinatorial technique they called \textit{reverse engineering}. There is an alternate $(4,3)$ factorization obtained by Hamada, which appears in \cite{nakamura} and is used to construct a factorization of type $(10,10)$. Since Xiao's fibration has 7 singular fibers, it must be of type $(4,3)$. The vanishing cycles of Xiao's fibration were determined recently by Szab\'{o} \cite{szabo}. 

Sato \cite{sato} showed that any genus-2 Lefschetz fibration with 7 singular fibers has total space diffeomorphic to $(S^2\times T^2)\#3\overline\CP$. Hence the fibrations of Baykur-Korkmaz, Hamada, and Xiao are supported on the same space. We show that the fibration structures are also the same:

\begin{thm}
\label{thm:equiv_factorizations}
The Lefschetz fibrations corresponding to the $(4,3)$ factorizations of Baykur-Korkmaz, Hamada, and Xiao are all isomorphic.
\end{thm}

We review the notion of isomorphic Lefschetz fibrations in Section \ref{sec:LF} and prove the theorem in Section \ref{sec:equiv}. Section \ref{sec:figs} contains most of our figures, including the vanishing cycles for the various factorizations.
\section{Lefschetz fibrations}
\label{sec:LF}
Let $X$ be a closed smooth oriented 4-manifold. A \textit{Lefschetz fibration} on $X$ is a smooth surjective map $f:X\to S^2$ such that:
\begin{enumerate}
\item $f$ has finitely many critical points $p_1,\dots,p_k$.
\item The local form of $f$ around a critical point is 
\[ (z_1,z_2)\mapsto w=z_1z_2,\]
where $(z_1,z_2)$ are local orientation-preserving complex coordinates around $p_i$ and $w$ is a local orientation-preserving complex coordinate around $f(p_i)$.
\end{enumerate}
By perturbing $f$ if necessary, we can always assume the critical points $p_i$ lie in distinct singular fibers. If $g$ is the genus of a regular fiber of $f$, we say that $f$ is a \textit{genus-$g$ Lefschetz fibration}.

\begin{dfn}
Let $f:X\to S^2$ and $f':X'\to S^2$ be Lefschetz fibrations with the same genus. We say $f$ and $f'$ are \textit{isomorphic} if there are orientation-preserving diffeomorphisms $h:S^2\to S^2$ and $H:X\to X'$ such that $f'H=hf$.
\end{dfn}

Let $f:X\to S^2$ be a genus-$g$ Lefschetz fibration. Fix a regular fiber $F=f^{-1}(q)$. Moving from $q$ to a critical value $f(p_i)$ in $S^2$ corresponds to traveling in $X$ from $F$ to the singular fiber $F_i$ containing $p_i$, and the local form of $f$ around $p_i$ implies that $F_i$ is obtained from $F$ by collapsing a curve $c_i\subseteq F$ to a point. $c_i$ is called a \textit{vanishing cycle}. Moreover, the local monodromy around $f(p_i)$ is the right-handed Dehn twist about $c_i$, denoted $t_{c_i}$. The global monodromy of $f$ around a disk containing all the critical values $f(p_1),\dots,f(p_k)$ is the product of these Dehn twists, and since $f$ extends over $S^2$, we see that
\[ t_{c_1}\dots t_{c_k}=1\]
in the mapping class group $\Mod_g$ of $\Sigma_g$. Such a relation is called a \textit{positive factorization}.

In the case $g=2$, a positive factorization $t_{c_1}\dots t_{c_k}=1$ in $\Mod_2$ with $n$ nonseparating curves and $s$ separating curves is said to be of \textit{type $(n,s)$}.

Conversely, a positive factorization $t_{c_1}\dots t_{c_k}=1$ in $\Mod_g$ determines a genus-$g$ Lefschetz fibration. Start with $\Sigma_g\times D^2$, attach $2$-handles along the curves $c_1,\dots,c_k$ with framing $-1$ with respect to the fiber framing, and attach another $\Sigma_g\times D^2$ to get a closed oriented 4-manifold $X$ equipped with a Lefschetz fibration $X\to S^2$. See Chapter 8 of \cite{gs} for details.

Obtaining a positive factorization from a Lefschetz fibration $f:X\to S^2$ requires some choices - a regular fiber $F=f^{-1}(q)$, an identification of $F$ with $\Sigma_g$, and paths in $S^2$ from $q$ to the critical values $f(p_i)$. It turns out \cite{mat} that when $g\ge 2$, there is a one-to-one correspondence between Lefschetz fibrations up to isomorphism and positive factorizations in $\Mod_g$, up to:
\begin{enumerate}
\item \textit{Hurwitz moves}: trading $t_{c_i}t_{c_{i+1}}$ with $t_{t_{c_i}(c_{i+1})}t_{c_i}$, and
\item \textit{global conjugations}: trading every $t_{c_i}$ for $t_{\phi(c_i)}$ for some $\phi\in\Mod_g$.
\end{enumerate}

\section{Computations in $S_{0,6}$}
Let $\iota:\Sigma_2\to\Sigma_2$ be the standard hyperelliptic involution on $\Sigma_2$ and let $p:\Sigma_2\to S_{0,6}$ be the resulting quotient map with $A$ the set of 6 branch points in $S_{0,6}$. We prefer to work with the images of curves in $S_{0,6}$ because they have fewer intersections, which simplifies Dehn twist computations. In this section, we recall some facts about curves in $\Sigma_2$ and their images in $S_{0,6}$, as well as how Dehn twists on $\Sigma_2$ correspond to maps on $S_{0,6}$.

We say a curve $\gamma\in\Sigma_2$ is \textit{hyperelliptic} if $\iota(\gamma)=\gamma$ setwise. Recall that every curve in $\Sigma_2$ is isotopic to a hyperelliptic curve.
\begin{lem}
Let $\gamma\subseteq\Sigma_2$ be a hyperelliptic curve.
\begin{enumerate}
\item If $\gamma$ is nonseparating, then $p(\gamma)$ is an arc with endpoints on $A$.
\item If $\gamma$ is separating, then $p(\gamma)$ is a simple closed curve disjoint from $A$.
\end{enumerate}
\end{lem}

If $\alpha$ is an arc in $S_{0,6}$ with endpoints on $A$, let $w_\alpha$ be the right half-twist about $\alpha$.
\begin{prop}
Let $c\subseteq\Sigma_2$ be a hyperelliptic curve, and let $\bar c=p(c)$ be its image in $S_{0,6}$.
\begin{enumerate}
\item If $c$ is nonseparating, then $t_c$ on $\Sigma_2$ corresponds to $w_{\bar c}$ on $S_{0,6}$. That is, $pt_c=w_{\bar c}p$.
\item If $c$ is separating, then $t_c$ on $\Sigma_2$ corresponds to $t_{\bar c}^2$ on $S_{0,6}$. That is, $pt_c=t_{\bar c}^2p$.
\end{enumerate}
\end{prop}

\begin{figure}[h]
\centering{
\resizebox{130mm}{!}{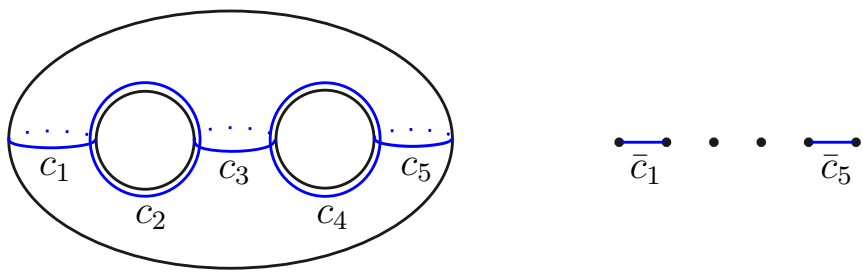}
\caption{The standard chain of curves on $\Sigma_2$.}
\label{fig:std_chain}
}
\end{figure}

Throughout, we will make use of the standard chain of curves on $\Sigma_2$, as shown in Figure \ref{fig:std_chain}.

\section{Equivalences}
\label{sec:equiv}
\begin{figure}[h]
\centering{
\resizebox{90mm}{!}{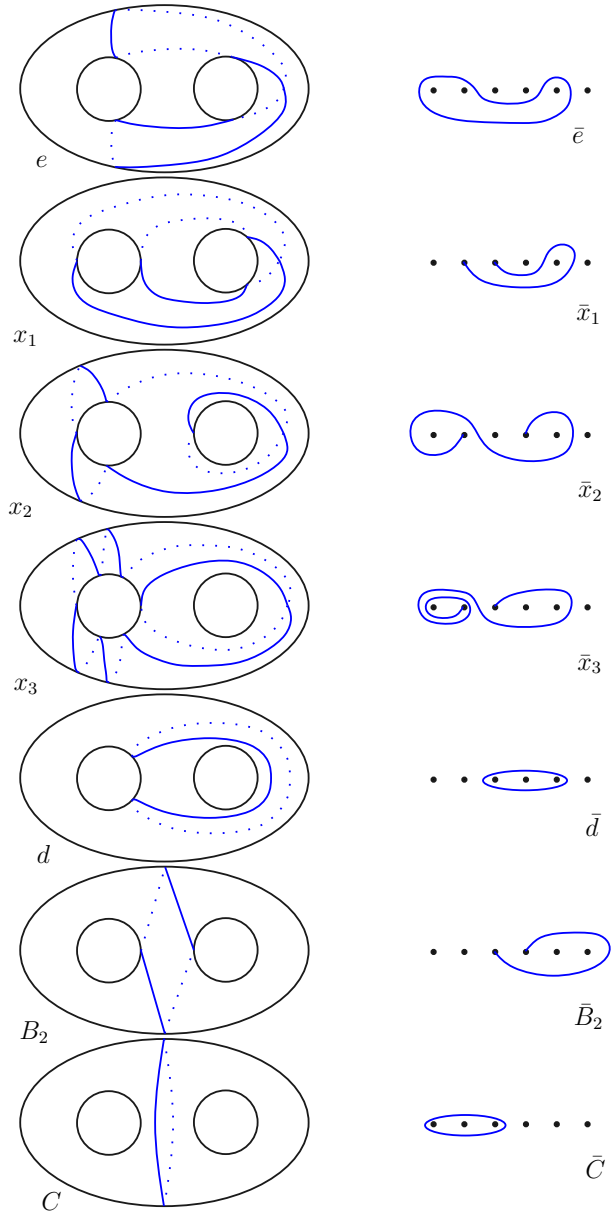}
\caption{Vanishing cycles in Baykur-Korkmaz's factorization.}
\label{fig:BK_cycles}
}
\end{figure}

In this section, we prove Theorem \ref{thm:equiv_factorizations}. The Baykur-Korkmaz $(4,3)$ factorization is
\begin{equation}
\label{eq:BK}
	t_et_{x_1}t_{x_2}t_{x_3}t_dt_{B_2}t_C=1,
\end{equation}
where the vanishing cycles and their images under $p$ appear in Figure \ref{fig:BK_cycles}. Likewise, the Hamada $(4,3)$ factorization is
\begin{equation}
\label{eq:H}
	t_\alpha t_D t_\sigma t_E t_\gamma t_F t_G=1,
\end{equation}
with the vanishing cycles and their images under $p$ in Figure \ref{fig:H_cycles}.

\begin{figure}[h]
\centering{
\resizebox{130mm}{!}{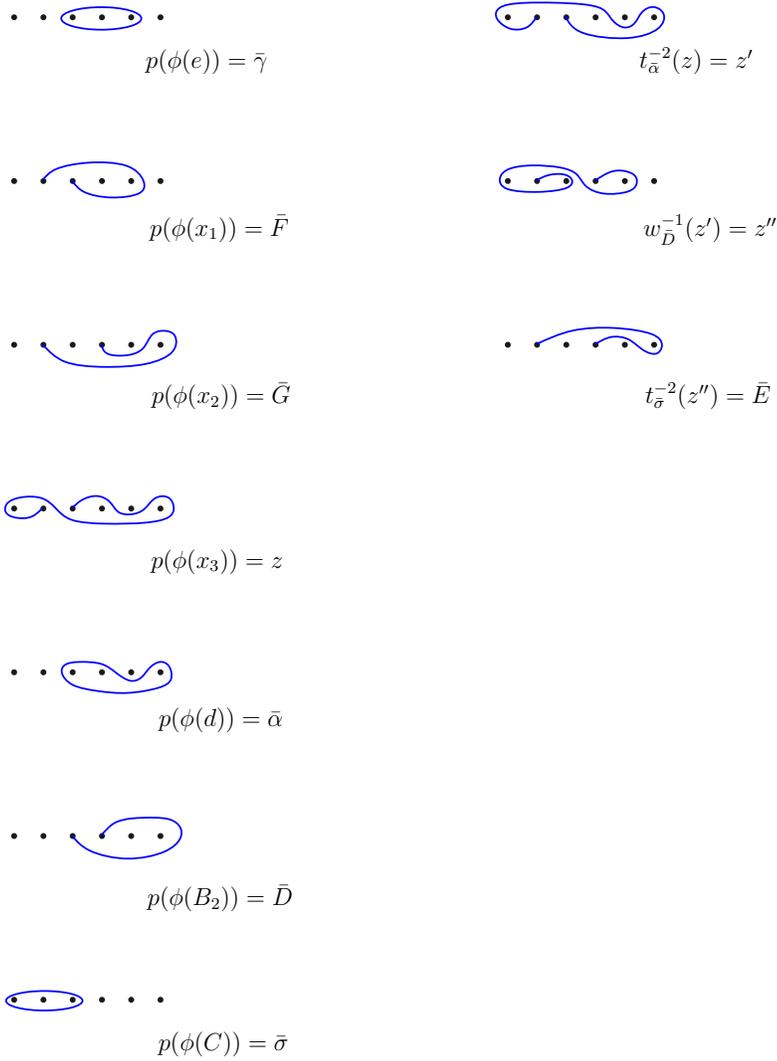}
\caption{Left: global conjugation of Baykur-Korkmaz by $\phi=t_{c_1}^{-2}t_{c_5}$. Right: Hurwitz moves to bring $t_{\phi(x_3)}$ to the right.}
\label{fig:BK_H_equiv}
}
\end{figure}

Let us start by showing the equivalence between (\ref{eq:BK}) and (\ref{eq:H}). Doing a global conjugation by $\phi=t_{c_1}^{-2}t_{c_5}$ to the vanishing cycles in (\ref{eq:BK}) results in the curves on the left of Figure \ref{fig:BK_H_equiv}. This changes the Baykur-Korkmaz factorization to the equivalent factorization
\[
	t_\gamma t_F t_G t_{\phi(x_3)}{\color{red}t_\alpha t_Dt_\sigma}=1.
\]
Moving the Dehn twists $t_\alpha, t_D$, and $t_\sigma$ (in red) to the left of $t_{\phi(x_3)}$ using Hurwitz moves turns $t_{\phi(x_3)}$ into $t_E$ (see the right of Figure \ref{fig:BK_H_equiv}). This results in the factorization
\[
	t_\gamma t_Ft_Gt_\alpha t_Dt_\sigma t_E=1,
\]
which after a cyclic permutation is precisely Hamada's factorization (\ref{eq:H}).

\begin{figure}[h]
\centering{
\resizebox{100mm}{!}{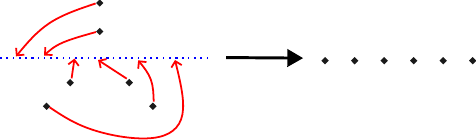}
\caption{Isotoping branch points to be standard.}
\label{fig:branch_isotope}
}
\end{figure}

\begin{figure}[h]
\centering{
\resizebox{100mm}{!}{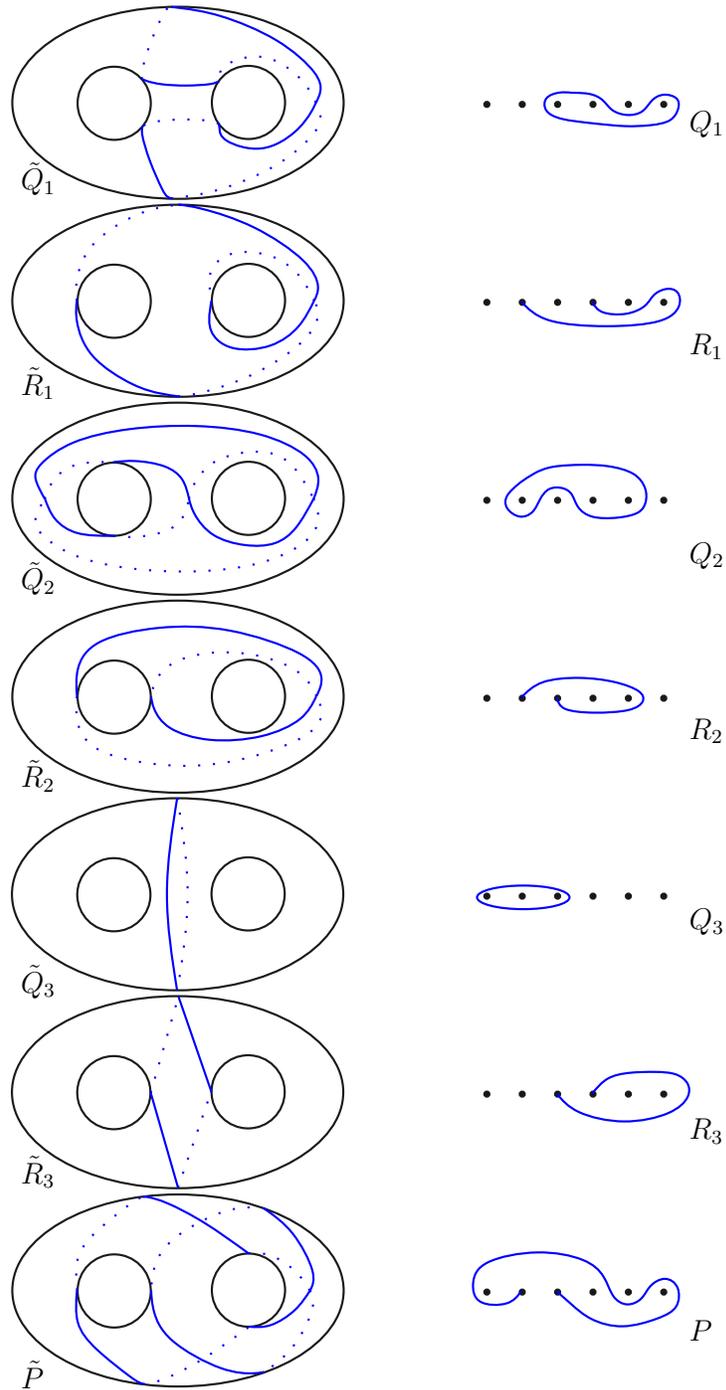}
\caption{Vanishing cycles in Xiao's $(4,3)$ fibration drawn in a standard way, with lifts to $\Sigma_2$.}
\label{fig:X_cycles}
}
\end{figure}

The vanishing cycles of Xiao's $(4,3)$ fibration were worked out by Szab\'{o} \cite{szabo}, and are given by their images under the quotient map $p:\Sigma_2\to S_{0,6}$. See Figure \ref{fig:X_cycles_planar}. We can isotope the branch points to be standard, as in Figure \ref{fig:branch_isotope}, and lift these vanishing cycles to $\Sigma_2$, as shown in Figure \ref{fig:X_cycles}. Xiao's factorization then reads:
\begin{equation}
\label{eq:X}
	t_{\tilde P}t_{\tilde R_3}t_{\tilde Q_3}t_{\tilde R_2}t_{\tilde Q_2}t_{\tilde R_1}t_{\tilde Q_1}=1.
\end{equation}

\begin{figure}[h]
\centering{
\resizebox{130mm}{!}{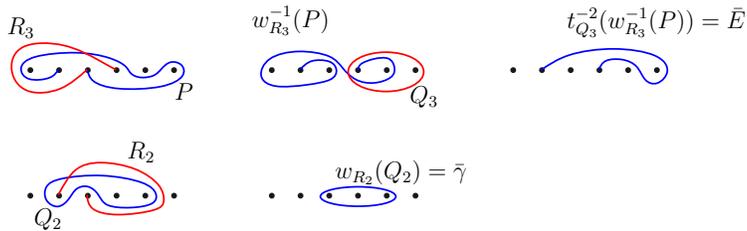}
\caption{Hurwitz moves for Xiao's factorization.}
\label{fig:X_equiv}
}
\end{figure}

We now show the equivalence between (\ref{eq:H}) and (\ref{eq:X}). First note that $\tilde R_3=D$, $\tilde Q_3=\sigma$, $\tilde R_2=F$, $\tilde R_1=G$, and $\tilde Q_1=\alpha$. Putting this into (\ref{eq:X}) and doing some Hurwitz moves (in red) gives:
\begin{equation*}
\begin{split}
	1&=t_{\tilde P}{\color{red}t_Dt_\sigma}t_F{\color{red}t_{\tilde Q_2}}t_Gt_\alpha\\
	&=t_Dt_\sigma t_Et_\gamma t_Ft_Gt_\alpha,
\end{split}
\end{equation*}
which, after a cyclic permutation, is precisely Hamada's factorization. The computations that $t_\sigma^{-1}t_D^{-1}(\tilde P)=E$ and $t_F(\tilde Q_2)=\gamma$ are shown in Figure \ref{fig:X_equiv}.

\begin{rmk}
\label{rmk:red_sections}
The factorizations of Baykur-Korkmaz, Hamada, and Xiao all have lifts to $\Mod(\Sigma_2^3)$ - the identity mapping class on the right sides of (\ref{eq:BK}), (\ref{eq:H}), and (\ref{eq:X}) are replaced with the product $t_{\delta_1}t_{\delta_2}t_{\delta_3}$ of boundary Dehn twists. Equivalently, the corresponding Lefschetz fibrations each have three sections of square $-1$. The boundaries are marked in red in Figure \ref{fig:sections}. There are many ways to see this - the Alexander method is probably the most straightforward. Alternatively, since each of the factorizations has hyperelliptic symmetry and each vanishing cycle misses the three red points, we automatically get lifts to $\Mod(\Sigma_2^3)$, following an argument in Section 4.2 of \cite{baykur}. 

There is yet another way to see this. In the process of obtaining the vanishing cycles of Xiao's fibration, Szab\'{o} located three sections \cite{szabo}, marked in red in Figure \ref{fig:X_cycles_planar}. After the isotopy in Figure \ref{fig:branch_isotope}, the sections end up exactly as in Figure \ref{fig:sections}.
\begin{figure}[h]
\centering{
\resizebox{110mm}{!}{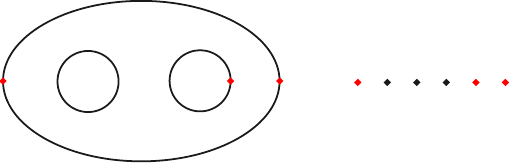}
\caption{Sections of $(4,3)$ factorizations.}
\label{fig:sections}
}
\end{figure}

The lifts of (\ref{eq:BK}), (\ref{eq:H}), and (\ref{eq:X}) to $\Mod(\Sigma_2^3)$ describe Lefschetz pencils with three base points (blow down the $(-1)$-sections). It is straightfoward to check that $t_\sigma^{-1}t_D^{-1}(\tilde P)=E$ and $t_F(\tilde Q_2)=\gamma$ holds in $\Sigma_2^3$, so our Hurwitz moves giving the equivalence between (\ref{eq:H}) and (\ref{eq:X}) also show that the corresponding Lefschetz pencils are equivalent.\footnote{The author would like to thank R {\.I}nan{\c{c}} Baykur for pointing this out.} Baykur-Hamada \cite{bh} have shown that the Lefschetz pencils corresponding to (\ref{eq:BK}) and (\ref{eq:H}) are equivalent.
\end{rmk}

\section{Reverse engineering Hamada's $(4,3)$ factorization}
It is possible, and not too tedious, to derive Hamada's $(4,3)$ factorization directly from the 4-chain relation using reverse engineering (see \cite{bk}). The idea is to take a standard relation with lots of nonseparating curves and find 2-chain substitutions. This cuts down the length of the relation by replacing 12 nonseparating curves with 1 separating curve.

We go in the opposite direction here - starting with Hamada's $(4,3)$ factorization, we replace separating Dehn twists with products of nonseparating Dehn twists, then use braid relations to arrive at a 4-chain relation. To simplify notation, we will write a curve to denote the Dehn twist along that curve, and write $t_i$ or even just $i$ to denote the Dehn twist along $c_i$ (as shown in Figure \ref{fig:std_chain}). For example, $\alpha 1E32$ really denotes $t_\alpha t_{c_1}t_Et_{c_3}t_{c_2}$.

Recall that Hamada's $(4,3)$ factorization is
\[ \alpha D\sigma E\gamma FG=I.\]
It is straightforward to check that
\begin{itemize}
\item $\alpha=t_4^{-1}t_3^{-1}(\sigma)$
\item $D=t_2t_1^2t_2(c_3)$
\item $E=t_4^2t_1^{-2}t_2^{-1}(c_3)$
\item $F=(t_3t_4)^{-3}(c_2)$
\item $G=t_4^{-2}t_1^2t_2(c_3).$
\end{itemize}
By the 2-chain relation, $\sigma=(12)^6$ and $\gamma=(34)^6$. It follows that
\begin{equation*}
\begin{split}
	I&=\alpha D{\color{red}(12)^6}E{\color{red}(34)^6}FG\\
	&=\alpha D(12)^3211{\color{red}211}E{\color{red}44}3443(34)^3FG\\
	&=\alpha D{\color{red}(12)^3}(211)(44)3(211)(3443){\color{red}(34)^3}FG\\
	&=\alpha D{\color{red}2112}11(211)(44)3(211)(3443)2{\color{red}(34)^3}G\\
	&=\alpha(2112)3(11)21144321134432(3443{\color{red}44})G\\
	&=\alpha 211231121144321134432(3443)G'4{\color{red}4}\\
	&=\alpha'42112311211443211344323443G'4,
\end{split}
\end{equation*}
where $G'=t_4^2(G)=t_1^2t_2(c_3)$ and $\alpha'=t_4(\alpha)=t_3^{-1}(\sigma)$. Continuing, we see that
\begin{equation*}
\begin{split}
	I&=\alpha'{\color{red}42112311}211443211344323443G'4\\
	&=\alpha'({\color{red}112}11243)211443211344323443G'4\\
	&=\alpha''(11243)211443211344323443G'4(112),
\end{split}
\end{equation*}
where $\alpha''=(t_1^2t_2)^{-1}(\alpha')=t_2^{-1}t_3^{-1}(\sigma)$. Observe that
\begin{equation*}
\begin{split}
	I&=\alpha''(11243)211443211344323443G'4({\color{red}112})\\
	&={\color{red}\alpha''}(11243)211443211344323443({\color{red}11}2)34\\
	&=1124321144321134432(11)344{\color{red}323}4\alpha''\\
	&=11243211443211344321134423{\color{red}2}4\alpha''\\
	&=211243211443211344321134{\color{red}4234}\alpha'\\
	&=211243211443211344321134(234{\color{red}3})\alpha'\\
	&=32112432114432113443211{\color{red}34234}(12)^6\\
	&=32112432114432113443211(23423)(12)^6\\
	&=321124321144321134W23(12)^6,
\end{split}
\end{equation*}
where $W=43211234$. Note that $W$ fixes $c_1,c_2$, and $c_3$. Hence $W$ commutes with $1,2$, and $3$. We are almost done:
\begin{equation*}
\begin{split}
	I&=321124321144321134W{\color{red}23}(12)^6\\
	&=3211243211{\color{red}44321134}(23)W(12)^6\\
	&=3211243211(23433211)23W(12)^6\\
	&=32112W{\color{red}3321123}W(12)^6\\
	&=321123321123WW(12)^6\\
	&=W'W'WW(12)^6,
\end{split}
\end{equation*}
where $W'=321123$. At the start of Section 3.2 in \cite{bk}, it is shown that $W'W'WW(12)^6=(1234)^{10}$. This demonstrates how to recover the 4-chain relation from Hamada's $(4,3)$ factorization, as desired.

\clearpage

\section{Vanishing cycles and Hurwitz moves}
\label{sec:figs}
\begin{figure}[h]
\centering{
\resizebox{90mm}{!}{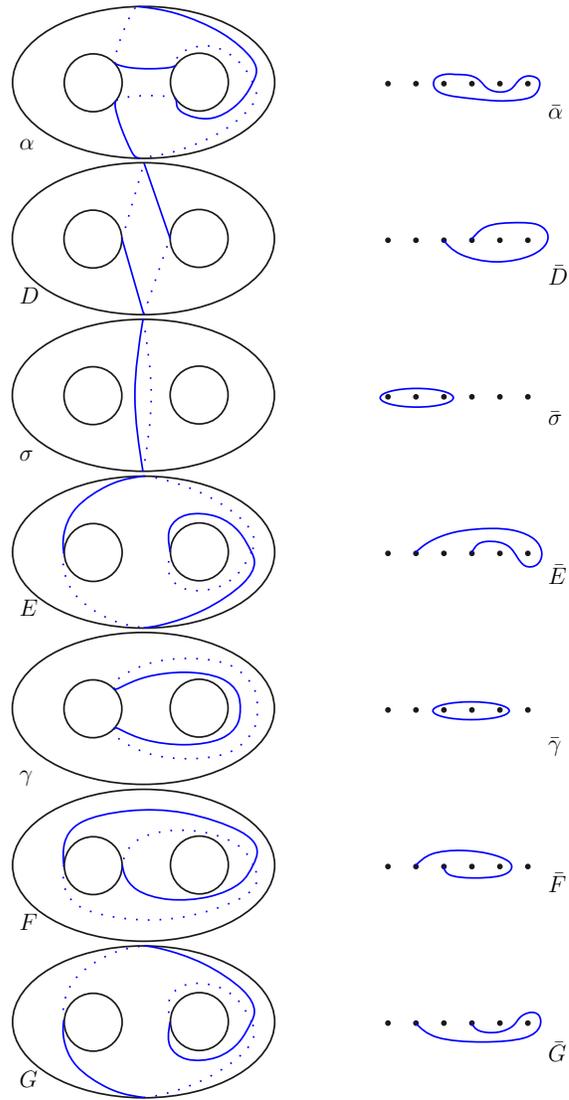}
\caption{Vanishing cycles in Hamada's factorization.}
\label{fig:H_cycles}
}
\end{figure}
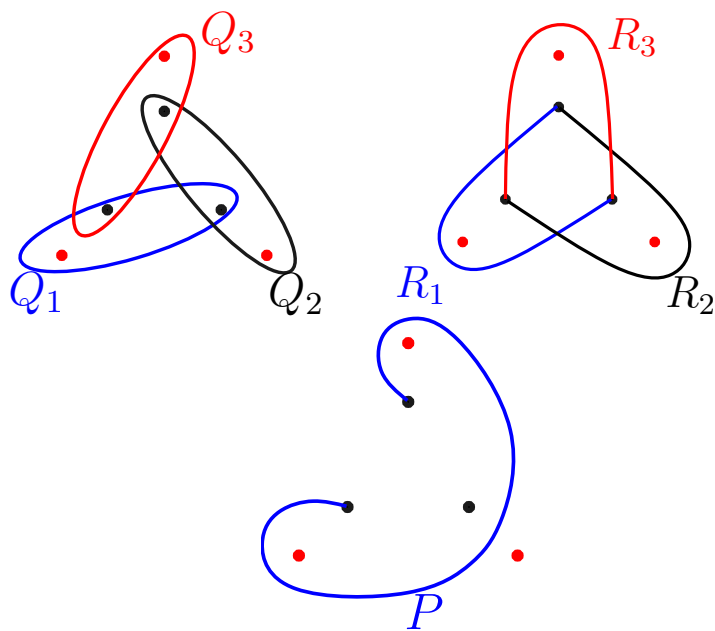
\begin{figure}[h]
\centering{
\resizebox{50mm}{!}{
\begingroup%
  \makeatletter%
  \providecommand\color[2][]{%
    \errmessage{(Inkscape) Color is used for the text in Inkscape, but the package 'color.sty' is not loaded}%
    \renewcommand\color[2][]{}%
  }%
  \providecommand\transparent[1]{%
    \errmessage{(Inkscape) Transparency is used (non-zero) for the text in Inkscape, but the package 'transparent.sty' is not loaded}%
    \renewcommand\transparent[1]{}%
  }%
  \providecommand\rotatebox[2]{#2}%
  \newcommand*\fsize{\dimexpr\f@size pt\relax}%
  \newcommand*\lineheight[1]{\fontsize{\fsize}{#1\fsize}\selectfont}%
  \ifx\svgwidth\undefined%
    \setlength{\unitlength}{93.63802091bp}%
    \ifx\svgscale\undefined%
      \relax%
    \else%
      \setlength{\unitlength}{\unitlength * \real{\svgscale}}%
    \fi%
  \else%
    \setlength{\unitlength}{\svgwidth}%
  \fi%
  \global\let\svgwidth\undefined%
  \global\let\svgscale\undefined%
  \makeatother%
  \begin{picture}(1,0.79103814)%
    \lineheight{1}%
    \setlength\tabcolsep{0pt}%
    \put(0,0){\includegraphics[width=\unitlength,page=1]{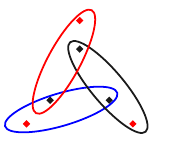}}%
    \put(-0.00559413,0.02488625){\color[rgb]{0,0,1}\makebox(0,0)[lt]{\lineheight{1.25}\smash{\begin{tabular}[t]{l}$Q_1$\end{tabular}}}}%
    \put(0.68491438,0.02012398){\color[rgb]{0,0,0}\makebox(0,0)[lt]{\lineheight{1.25}\smash{\begin{tabular}[t]{l}$Q_2$\end{tabular}}}}%
    \put(0.50395335,0.71539464){\color[rgb]{1,0,0}\makebox(0,0)[lt]{\lineheight{1.25}\smash{\begin{tabular}[t]{l}$Q_3$\end{tabular}}}}%
  \end{picture}%
\endgroup%
}
\resizebox{50mm}{!}{
\begingroup%
  \makeatletter%
  \providecommand\color[2][]{%
    \errmessage{(Inkscape) Color is used for the text in Inkscape, but the package 'color.sty' is not loaded}%
    \renewcommand\color[2][]{}%
  }%
  \providecommand\transparent[1]{%
    \errmessage{(Inkscape) Transparency is used (non-zero) for the text in Inkscape, but the package 'transparent.sty' is not loaded}%
    \renewcommand\transparent[1]{}%
  }%
  \providecommand\rotatebox[2]{#2}%
  \newcommand*\fsize{\dimexpr\f@size pt\relax}%
  \newcommand*\lineheight[1]{\fontsize{\fsize}{#1\fsize}\selectfont}%
  \ifx\svgwidth\undefined%
    \setlength{\unitlength}{99.94651975bp}%
    \ifx\svgscale\undefined%
      \relax%
    \else%
      \setlength{\unitlength}{\unitlength * \real{\svgscale}}%
    \fi%
  \else%
    \setlength{\unitlength}{\svgwidth}%
  \fi%
  \global\let\svgwidth\undefined%
  \global\let\svgscale\undefined%
  \makeatother%
  \begin{picture}(1,0.77539416)%
    \lineheight{1}%
    \setlength\tabcolsep{0pt}%
    \put(0,0){\includegraphics[width=\unitlength,page=1]{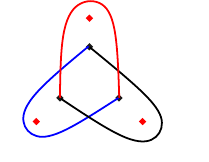}}%
    \put(-0.00524103,0.04421618){\color[rgb]{0,0,1}\makebox(0,0)[lt]{\lineheight{1.25}\smash{\begin{tabular}[t]{l}$R_1$\end{tabular}}}}%
    \put(0.7130682,0.01744691){\color[rgb]{0,0,0}\makebox(0,0)[lt]{\lineheight{1.25}\smash{\begin{tabular}[t]{l}$R_2$\end{tabular}}}}%
    \put(0.55691395,0.70452523){\color[rgb]{1,0,0}\makebox(0,0)[lt]{\lineheight{1.25}\smash{\begin{tabular}[t]{l}$R_3$\end{tabular}}}}%
  \end{picture}%
\endgroup%
}
\resizebox{35mm}{!}{
\begingroup%
  \makeatletter%
  \providecommand\color[2][]{%
    \errmessage{(Inkscape) Color is used for the text in Inkscape, but the package 'color.sty' is not loaded}%
    \renewcommand\color[2][]{}%
  }%
  \providecommand\transparent[1]{%
    \errmessage{(Inkscape) Transparency is used (non-zero) for the text in Inkscape, but the package 'transparent.sty' is not loaded}%
    \renewcommand\transparent[1]{}%
  }%
  \providecommand\rotatebox[2]{#2}%
  \newcommand*\fsize{\dimexpr\f@size pt\relax}%
  \newcommand*\lineheight[1]{\fontsize{\fsize}{#1\fsize}\selectfont}%
  \ifx\svgwidth\undefined%
    \setlength{\unitlength}{61.43932072bp}%
    \ifx\svgscale\undefined%
      \relax%
    \else%
      \setlength{\unitlength}{\unitlength * \real{\svgscale}}%
    \fi%
  \else%
    \setlength{\unitlength}{\svgwidth}%
  \fi%
  \global\let\svgwidth\undefined%
  \global\let\svgscale\undefined%
  \makeatother%
  \begin{picture}(1,1.21605142)%
    \lineheight{1}%
    \setlength\tabcolsep{0pt}%
    \put(0,0){\includegraphics[width=\unitlength,page=1]{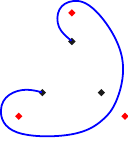}}%
    \put(0.54720584,0.0283817){\color[rgb]{0,0,1}\makebox(0,0)[lt]{\lineheight{1.25}\smash{\begin{tabular}[t]{l}$P$\end{tabular}}}}%
  \end{picture}%
\endgroup%
}
\caption{Vanishing cycles of Xiao's $(4,3)$ fibration. Sections are in red, see Remark \ref{rmk:red_sections}.}
\label{fig:X_cycles_planar}
}
\end{figure}


\end{document}